\documentclass[psamsfonts]{conm-p-l}
\numberwithin{equation}{section}

\def\lBrack{\mathopen{[\![}}
\def\rBrack{\mathclose{]\!]}}
\def\bB{{\mathbb B}}
\def\MBU{{M^\bB\!/U}}
\def\forces{\mathrel{||\hspace{-2.2pt}-}}
\def\OR{\mathrel{\mathrm{OR}}}
\def\AND{\mathrel{\mathrm{AND}}}
\def\bval#1{\lBrack #1 \rBrack^\bB}

\begin{document}

\title[A beginner's guide to forcing]{A beginner's guide to forcing}
\author{Timothy Y. Chow}
\address{Center for Communications Research, 805 Bunn Drive, Princeton, NJ 08540}
\email{tchow@alum.mit.edu}
\urladdr{http://alum.mit.edu/www/tchow}
\dedicatory{Dedicated to Joseph Gallian on his 65th birthday}
\date{December 14, 2007}

\maketitle

\section{Introduction}

In 1963, Paul Cohen stunned the mathematical world with his new
technique of \textit{forcing,} which allowed him to solve several
outstanding problems in set theory at a single stroke.
Perhaps most notably, he proved
the independence of the continuum hypothesis (CH) from the
Zermelo-Fraenkel-Choice (ZFC) axioms of set theory.
The impact of Cohen's ideas on the practice of set theory,
as well as on the philosophy of mathematics,
has been incalculable.

Curiously, though, despite the importance of Cohen's work
and the passage of nearly fifty years,
forcing remains totally mysterious to the vast majority of mathematicians,
even those who know a little mathematical logic.
As an illustration, let us note that
Monastyrsky's outstanding book~\cite{monastyrsky}
gives highly informative and insightful expositions
of the work of almost every Fields Medalist---but
says almost nothing about forcing.
Although there exist numerous textbooks with
mathematically correct and complete proofs of the
basic theorems of forcing,
the subject remains notoriously difficult for beginners to learn.

All mathematicians are familiar with the concept of an
\textit{open research problem.}
I propose the less familiar concept of an \textit{open exposition problem.}
Solving an open exposition problem means
explaining a mathematical subject in a way that 
renders it totally perspicuous.
Every step should be motivated and clear;
ideally, students should feel that they could have
arrived at the results themselves.
The proofs should be ``natural'' in Donald Newman's sense~\cite{newman}:
\begin{quotation}
This term $\ldots$ is introduced to mean not having any ad hoc
constructions or \textit{brilliancies.}  A ``natural'' proof, then,
is one which proves itself, one available to the ``common
mathematician in the streets.''
\end{quotation}

I believe that it is an open exposition problem
to explain forcing.
Current treatments allow readers to verify the truth of the basic theorems,
and to progress fairly rapidly to the point where they can \textit{use}
forcing to prove their own independence results
(see \cite{baumgartner} for a particularly nice explanation of
how to use forcing as a black box to turn independence questions
into concrete combinatorial problems).
However, in all treatments that I know of,
one is left feeling that only a genius with fantastic intuition
or technical virtuosity could have found the road to the final result.

This paper does not solve this open exposition problem,
but I believe it is a step in the right direction.
My goal is to give a rapid overview of the subject,
emphasizing the broad outlines and
the intuitive motivation while omitting most of the proofs.
The reader will not, of course, master forcing by reading this
paper in isolation without consulting standard textbooks for the omitted details,
but my hope is to provide a map of the forest
so that the beginner will not get lost while forging through the trees.
Currently, no such bird's-eye overview seems to be available in
the published literature; I hope to fill this gap.
I also hope that this paper will inspire others to continue the job of
making forcing totally transparent.

\section{Executive summary}

The negation of~CH says that there is a cardinal
number, $\aleph_1$, between the cardinal numbers $\aleph_0$ and
$2^{\aleph_0}$.  One might therefore try to build a structure
that satisfies the negation of~CH by starting with something
that \textit{does} satisfy CH
(G\"odel had in fact constructed such structures)
and ``inserting'' some sets that are missing.

The fundamental theorem of forcing is that,
under very general conditions, one can indeed start with
a mathematical structure~$M$ that satisfies the ZFC axioms,
and enlarge it by adjoining a new element~$U$
to obtain a new structure~$M[U]$ that also satisfies ZFC.
Conceptually, this process is analogous to the process of adjoining
a new element~$X$ to, say, a given ring~$R$
to obtain a larger ring~$R[X]$.
However, the construction of~$M[U]$ is a lot more complicated
because the axioms of ZFC are more complicated than the axioms for a ring.
Cohen's idea was to build the new element~$U$ one step at a time,
tracking what new properties of~$M[U]$ would be ``forced'' to hold
at each step, so that one could control the properties
of~$M[U]$---in particular, making it satisfy the negation of~CH as well as
the axioms of ZFC.

The rest of this paper fleshes out the above construction in more detail.

\section{Models of ZFC}

As mentioned above, Cohen proved the independence of CH from ZFC;
more precisely, he proved that if ZFC is consistent,
then CH is not a logical consequence of the ZFC axioms.
G\"odel had already proved that if ZFC is consistent,
then $\neg$CH, the negation of CH,
is not a logical consequence of ZFC,
using his concept of ``constructible sets.''
(Note that the hypothesis that ZFC is consistent cannot be dropped,
because if ZFC is inconsistent then
\textit{everything} is a logical consequence of~ZFC!)

Just how does one go about proving that CH is not a logical consequence of~ZFC?
At a very high level,
the structure of the proof is what you would expect:
One writes down a very precise statement of the ZFC axioms
and of $\neg$CH,
and then one constructs a mathematical structure that satisfies both ZFC
and $\neg$CH.
This structure is said to be a \textit{model} of the axioms.
Although
the term ``model'' is not often seen in mathematics outside of formal logic,
it is actually a familiar concept.
For example, in group theory, a ``model of the group-theoretic axioms''
is just a group, i.e., a set $G$ with a binary operation $*$
satisfying axioms such as: ``There exists an element $e$ in~$G$
such that $x * e = e * x = x$ for all $x$ in $G$,'' and so forth.

Analogously, we could invent a term---say, \textit{universe}---to mean
``a structure that is a model of ZFC.''
Then we could begin our study of ZFC with definition such as,
``A \textit{universe} is a set~$M$ together with a binary relation~$R$
satisfying$\dots$'' followed by a long list of axioms such as
the \textit{axiom of extensionality:}

\begin{quotation}
If $x$ and $y$ are distinct elements of~$M$ then either there
exists $z$ in $M$ such that $z R x$ but not $z R y$, or there
exists $z$ in $M$ such that $z R y$ but not $z R x$.
\end{quotation}
Another axiom of ZFC is the \textit{powerset axiom:}

\begin{quotation}
For every $x$ in~$M$,
there exists $y$ in~$M$ with the following property:
For every $z$ in~$M$, $z R y$ if and only if $z\subseteq x$.
\end{quotation}
(Here the expression ``$z\subseteq x$''
is an abbreviation for ``every $w$ in~$M$ satisfying $w R z$
also satisfies $w R x$.'')
There are other axioms, which can be found in any set theory textbook,
but the general idea should be clear from these two examples.
Note that the binary relation is usually denoted by the symbol~$\in$
since the axioms are inspired by the set membership relation.
However, we have deliberately chosen the unfamiliar symbol~$R$
to ensure that the reader will not misinterpret the axiom
by accidentally reading $\in$ as ``is a member of.''

As an aside, we should mention that it is not standard
to use the term \textit{universe} to mean ``model of ZFC.''
For some reason set theorists tend to give
a snappy name like ``ZFC'' to a \textit{list of axioms,}
and then use the term ``model of ZFC'' to refer to
the \textit{structures that satisfy the axioms},
whereas in the rest of mathematics it is the other way around:
one gives a snappy name like ``group'' to the structure, and 
then uses the term ``axioms for a group'' to refer to the axioms.
Apart from this terminological point, though, the formal setup here
is entirely analogous to that of group theory.
For example, in group theory, the statement~$S$ that
``$x*y = y*x$ for all $x$ and~$y$''
is not a logical consequence of the axioms of group theory,
because there exists a mathematical structure---namely
a non-abelian group---that satisfies the group axioms
as well as the negation of~$S$.

On the other hand, the definition of a model of ZFC has some curious
features, so a few additional remarks are in order.

\subsection{Apparent circularity}

One common confusion about models of ZFC
stems from a tacit expectation that some people have,
namely that we are supposed to suspend all our preconceptions
about sets when beginning the study of ZFC.
For example, it may have been a surprise to some readers to see
that a universe is defined to be a \textit{set} together with$\ldots\,$.
Wait a minute---what is a set?
Isn't it circular to define sets in terms of sets?

In fact, we are not defining sets in terms of sets,
but \textit{universes} in terms of sets.
Once we see that all we are doing is studying a subject
called ``universe theory'' (rather than ``set theory''),
the apparent circularity disappears.

The reader may still be bothered by the lingering feeling
that the point of introducing ZFC is to ``make set theory rigorous''
or to examine the foundations of mathematics.
While it is true that ZFC can be used as a tool
for such philosophical investigations,
we do not do so in this paper.
Instead, we take for granted that ordinary mathematical
reasoning---including reasoning about sets---is perfectly valid and
does not suddenly become invalid when the object of study is ZFC.
That is, we approach the study of ZFC and its models
in the same way that one approaches the study of any other
mathematical subject.
This is the best way to grasp the mathematical content;
after this is achieved, one can then try to apply the technical
results to philosophical questions if one is so inclined.

Note that in accordance with our attitude that
ordinary mathematical reasoning is perfectly valid,
we will freely employ reasoning about infinite sets
of the kind that is routinely used in mathematics.
We reassure readers who harbor philosophical doubts about the validity of
infinitary set-theoretic reasoning 
that Cohen's proof can be turned into a purely finitistic one.
We will not delve into such metamathematical niceties here,
but see for example the beginning of Chapter VII of Kunen's book~\cite{kunen}.

\subsection{Existence of examples}

A course in group theory typically begins with many examples of groups.
One then verifies that the examples satisfy all the axioms of group theory.
Here we encounter an awkward feature of models of ZFC,
which is that exhibiting explicit models of ZFC is difficult.
For example, there are no finite models of ZFC.
Worse, by a result known as the \textit{completeness theorem,}
the statement that ZFC has \textit{any models at all} is
equivalent to the statement that ZFC is consistent,
which is an assumption that is at least mildly controversial.
So how can we even get off the ground?

Fortunately, these difficulties are not as severe as they might seem at first.
For example, one entity that is almost a model of ZFC is~$V$,
the class of all sets.
If we take $M=V$ and we take $R$ to mean ``is a member of,''
then we see that the axiom of extensionality simply says that
two sets are equal if and only if
they contain the same elements---a manifestly true statement.
The rest of the axioms of ZFC are similarly self-evident
when $M=V$.\footnote{Except, perhaps, the axiom of regularity,
but this is a technical quibble that we shall ignore.}
The catch is that a model of ZFC has to be a \textit{set,}
and $V$, being ``too large'' to be a set (Cantor's paradox),
is a proper class and therefore,
strictly speaking, is disqualified from being a model of ZFC.
However, it is close enough to being a model of ZFC to be
intuitively helpful.

As for possible controversy over whether ZFC is consistent,
we can sidestep the issue simply by treating the consistency of ZFC
like any other unproved statement, such as the Riemann hypothesis.
That is, we can assume it freely as long as we remember to preface
all our theorems with a conditional clause.\footnote{In fact,
we already did this when we said that the precise statement of
Cohen's result is that \textit{if ZFC is consistent} then
CH is not a logical consequence of ZFC.}
So from now on we shall assume that ZFC is consistent,
and therefore that models of ZFC exist.

\subsection{``Standard'' models}

Even granting the consistency of ZFC,
it is not easy to produce models.
One can extract an example from the
proof of the completeness theorem, but this example is unnatural
and is not of much use for tackling~CH.
Instead of continuing the search for explicit examples,
we shall turn our attention to important properties of models of~ZFC.

One important insight of Cohen's was that it is useful to consider
what he called \textit{standard} models of ZFC.
A model $M$ of ZFC is standard if the elements of~$M$
are \textit{well-founded sets} and if the relation $R$ is
ordinary set membership.
Well-founded sets are sets that are built up inductively from
the empty set, using operations such as taking unions, subsets, powersets, etc.
Thus the empty set $\{\}$ is well-founded, as are
$\{\{\}\}$ and the infinite set $\{\{\}, \{\{\}\}, \{\{\{\}\}\}, \ldots\}$.
They are called ``well-founded'' because the nature of their
inductive construction precludes any well-founded set from being a
member of itself.
We emphasize that if $M$ is standard, then the elements of~$M$ are not
amorphous ``atoms,'' as some of us envisage the elements of an
abstract group to be, but are~\textit{sets.}
Moreover, well-founded sets are \textit{not}
themselves built up from ``atoms''; it's ``sets all the way down.''

While it is fairly clear that if standard models of ZFC exist,
then they form a natural class of examples,
it is not at all clear that any standard models exist at all,
even if ZFC is consistent.\footnote{It turns out that
the existence of a standard model of ZFC is indeed a stronger
assumption than the consistency of ZFC, but we will ignore this nicety.}
(The class of \text{all} well-founded sets
is a proper class and not a set and hence is disqualified.)
Moreover, even if standard models exist, one might think that
constructing a model of ZFC satisfying $\neg$CH might require
considering ``exotic'' models in which the binary relation~$R$
bears very little resemblance to ordinary set membership.
Cohen himself admits on page~108 of~\cite{cohen} that
a minor leap of faith is involved here:

\begin{quotation}
Since the negation of CH or AC may appear to be somewhat unnatural
one might think it hopeless to look for standard models.
However, we make a firm decision at the point to consider only
standard models.  Although this may seem like a very severe
limitation in our approach it will turn out that this very limitation
will guide us in suggesting possibilities.
\end{quotation}

Another property that a model of ZFC can have is
\textit{transitivity.}  A standard model $M$ of ZFC is \textit{transitive}
if every member of an element of~$M$ is also an element of~$M$.
(The term \textit{transitive} is used because we can write
the condition in the suggestive form ``$x\in y$ and $y\in M$
implies $x\in M$.'')
This is a natural condition to impose if we think of $M$ as a
universe consisting of ``all that there is'';
in such a universe, sets ``should'' be sets of things
that already exist in the universe.
Cohen's remark about standard models
applies equally to transitive models.\footnote{We remark in passing
that the \textit{Mostowski collapsing theorem} implies that
if there exist \textit{any} standard models of ZFC,
then there exist standard transitive models.}

Our focus will be primarily on standard transitive models.
Of course, this choice of focus is made with the benefit of hindsight,
but even without the benefit of hindsight,
it makes sense to study models with natural properties before
studying exotic models.
When $M$ is a standard transitive model, we will often use the
symbol $\in$ for the relation~$R$, because in this case $R$
\textit{is} in fact set membership.

\section{Powersets and absoluteness}

At some point in their education, most mathematicians learn
that all familiar mathematical objects can be defined in terms of sets.
For example, one can define the number~$0$ to be the empty set~$\{\}$,
the number~$1$ to be the set~$\{0\}$, and in general the number~$n$
to be the so-called \textit{von Neumann ordinal} $\{0, 1, \ldots, n-1\}$.
The set~$\mathbb N$ of all natural numbers may be defined to be~$\aleph_0$,
the set of all von Neumann ordinals.\footnote{We
elide the distinction between the cardinality~$\aleph_0$ of~$\mathbb N$
and the order type~$\omega$ of~$\mathbb N$.}
Note that with these definitions, the membership relation on~$\aleph_0$
corresponds to the usual ordering on the natural numbers
(this is why $n$ is defined as $\{0,1,\ldots,n-1\}$ rather
than as~$\{n-1\}$).
The ordered pair $(x,y)$ may be defined \`a~la Kuratowski as
the set $\{\{x\}, \{x,y\}\}$.
Functions, relations, bijections, maps, etc., can be defined as
certain sets of ordered pairs.
More interesting mathematical structures can be defined as
ordered pairs $(X,S)$ where $X$ is an underlying set and
$S$ is the structure on~$X$.
With this understanding, the class~$V$ of all sets
may be thought of as being the entire mathematical universe.

Models of~ZFC, like everything else, live inside~$V$,
but they are special because they look a lot like~$V$ itself.
This is because it turns out that virtually all mathematical proofs
of the existence of some object~$X$ in~$V$
can be mimicked by a proof from the ZFC axioms,
thereby proving that any model of ZFC must contain an object
that is at least highly analogous to~$X$.
It turns out that this ``analogue'' of~$X$ is often \textit{equal} to~$X$,
especially when $M$ is a standard transitive model.
For example, it turns out that every standard transitive model $M$ of ZFC
contains all the von Neumann ordinals as well as~$\aleph_0$.

However, the analogue of a mathematical object~$X$ is not \textit{always}
equal to~$X$.
A crucial counterexample is the powerset of~$\aleph_0$,
denoted by~$2^{\aleph_0}$.
Na\"\i vely, one might suppose that the powerset axiom of ZFC
guarantees that $2^{\aleph_0}$ must be a member of
any standard transitive model~$M$.
But let us look more closely at the precise statement
of the powerset axiom.
Given that $\aleph_0$ is in~$M$,
the powerset axiom guarantees the existence of $y$ in~$M$
with the following property:
For every $z$ in~$M$, $z\in y$ if and only if
every $w$ in~$M$ satisfying $w \in z$
also satisfies $w \in \aleph_0$.
Now, does it follow that $y$ is precisely
the set of all subsets of~$\aleph_0$?

No.  First of all, it is not even immediately clear
that $z$ is a subset of~$\aleph_0$;
the axiom does not require that \textit{every $w$}
satisfying $w\in z$ also satisfies $w\in \aleph_0$;
it requires only that \textit{every $w$ in~$M$}
satisfying $w\in z$ satisfies $w\in x$.
However, under our assumption that $M$ is \textit{transitive,}
every $w \in z$ is in fact in~$M$, so indeed $z$ is a subset of~$\aleph_0$.

More importantly, though, $y$ does not contain \textit{every}
subset of~$\aleph_0$; it contains
\textit{only those subsets of $x$ that are in~$M$.}
So if, for example, $M$ happens to be \textit{countable}
(i.e., $M$ contains only countably many elements),
then $y$ will be countable, and so a fortiori
$y$ cannot be equal to $2^{\aleph_0}$,
since $2^{\aleph_0}$ is uncountable.
The set $y$, which we might call the \textit{powerset of $\aleph_0$ in~$M$,}
is not the same as the ``real''
powerset of~$\aleph_0$, a.k.a.~$2^{\aleph_0}$;
many subsets of~$\aleph_0$ are ``missing'' from~$y$.

This is a subtle and important point, so let us explore it further.
We may ask, is it really possible for a standard transitive model of ZFC
to be countable?
Can we not mimic (in ZFC) Cantor's famous proof that
$2^{\aleph_0}$ is uncountable to show that
$M$ must contain an uncountable set,
and then conclude by transitivity that $M$ itself must be uncountable?

The answer is no.  Cantor's theorem states that
there is no bijection between $\aleph_0$ and $2^{\aleph_0}$.
If we carefully mimic Cantor's proof with a proof from the ZFC axioms,
then we find that Cantor's theorem tells us
that there is indeed a set $y$ in~$M$ that plays the role
of the powerset of~$\aleph_0$ in~$M$, and that
there is \textit{no bijection in~$M$} between $\aleph_0$ and~$y$.
However, this fact does not mean that there is \textit{no bijection at all}
between $\aleph_0$ and~$y$.
There might be a bijection \textit{in~$V$} between them;
we know only that such a bijection \textit{cannot be a member of~$M$};
it is ``missing'' from~$M$.
So Cantor's theorem does not exclude the possibility that~$y$,
as well as~$M$, is countable, even though $y$ is necessarily
``uncountable in~$M$.''\footnote{This curious state of affairs
often goes by the name of \textit{Skolem's paradox.}}
It turns out that something stronger can be said:
the so-called \textit{L\"owenheim-Skolem theorem} says that
if there are any models of ZFC at all, then in fact there
exist countable models.

More generally, one says that a concept in~$V$ is \textit{absolute}
if it coincides with its counterpart in~$M$.
For example, ``the empty set,'' ``is a member of,''
``is a subset of,'' ``is a bijection,''
and ``$\aleph_0$'' all turn out to be absolute
for standard transitive models.
On the other hand, ``is the powerset of''
and ``uncountable'' are \textit{not} absolute.
For a concept that is not absolute,
we must distinguish carefully between the concept ``in the real world''
(i.e., in~$V$) and the concept in~$M$.

A careful study of ZFC necessarily requires keeping track
of exactly which concepts are absolute and which are not.
However, since the majority of basic concepts are absolute,
except for those associated with taking powersets and cardinalities,
in this paper we will adopt the approach of mentioning non-absoluteness
only when it is especially relevant.

\section{How one might try to build a model satisfying $\neg$CH}
\label{sec:CH}

The somewhat counterintuitive fact that ZFC has countable models
with many missing subsets provides a hint as to
how one might go about constructing
a model for ZFC that satisfies $\neg$CH.
Start with a countable standard transitive model~$M$.
The elementary theory of cardinal numbers tells us that
there is always a smallest cardinal number after any given cardinal number,
so let $\aleph_1, \aleph_2, \ldots$ denote the next largest cardinals
after $\aleph_0$.
As usual we can mimic the proofs of these facts about cardinal numbers
with formal proofs from the axioms of~ZFC,
to conclude that there is a set in~$M$ that plays the role of~$\aleph_2$ in~$M$.
We denote this set by $\aleph_2^M$.
Let us now construct a function $F$ from
the Cartesian product $\aleph_2^M\times \aleph_0$ into the set $2 = \{0,1\}$.
We may interpret $F$ as a sequence of functions from $\aleph_0$ into~$2$.
Because $M$ is countable and transitive, so is $\aleph_2^M$;
thus we can easily arrange for these functions to be pairwise distinct.
Now, if $F$ is already in~$M$, then $M$ satisfies $\neg$CH!
The reason is that functions from $\aleph_0$ into~$2$
can be identified with subsets of~$\aleph_0$,
and $F$ therefore shows us that the powerset of~$\aleph_0$ in~$M$
must be at least $\aleph_2$ in~$M$.  Done!

But what if $F$ is missing from~$M$?
A natural idea is to add $F$ to~$M$ to obtain a larger model of ZFC,
that we might call $M[F]$.\footnote{Later on we will use the notation
$M[U]$ rather than $M[F]$ because it will turn out to be more
natural to think of the larger model as being obtained by adjoining
another set~$U$ that is closely related to~$F$,
rather than by adjoining $F$ itself.
For our purposes, $M[U]$ and $M[F]$ can just be thought of as two
different names for the same object.}
The hope would be that $F$ can be added in a way that does not
``disturb'' the structure of~$M$ too much, so that the argument
in the previous paragraph can be carried over into $M[F]$,
which would therefore satisfy~$\neg$CH.

Miraculously, this seemingly na\"\i ve idea actually works!
There are, of course, numerous technical obstacles to be surmounted,
but the basic plan as outlined above is on the right track.
For those who like to think algebraically,
it is quite appealing to learn that forcing is a technique
for constructing new models from old ones by adjoining a new element
that is missing from the original model.
Even without any further details,
one can already imagine that the ability to adjoin new elements
to an existing model
gives us enormous flexibility in our quest to create models with
desired properties.
And indeed, this is true; it is the reason why forcing is such a
powerful idea.

What technical obstacles need to be surmounted?
The first thing to note is that one clearly cannot
add only the set $F$ to~$M$ and expect to obtain a model of ZFC;
one must also add, at minimum,
every set that is ``constructible'' from~$F$ together with elements of~$M$,
just as when we create an extension of an algebraic object by
adjoining~$x$, we must also adjoin everything that is
\textit{generated} by~$x$.
We will not define ``constructible'' precisely here,
but it is the same concept
that G\"odel used to prove that CH is consistent with ZFC,
and in particular it was already a familiar concept before
Cohen came onto the scene.

A more serious obstacle is that it turns out that we cannot,
for example, simply take an arbitrary subset~$a$ of~$\aleph_0$
that is missing from~$M$ and adjoin~$a$, along with everything
constructible from~$a$ together with elements of~$M$, to~$M$;
the result will not necessarily be a model of ZFC.
A full explanation of this result would take us
too far afield---the interested reader should see
page~111 of Cohen's book~\cite{cohen}---but
the rough idea is that we could perversely choose $a$ to
be a set that encodes explicit information about the size of~$M$,
so that adjoining~$a$ would create a kind of self-referential paradox.
Cohen goes on to say:

\begin{quotation}
Thus $a$ must have certain special properties$\ldots\,$.
Rather than describe $a$ directly,
it is better to examine the various properties of~$a$
and determine which are desirable and which are not.
The chief point is that we do not wish $a$ to contain
``special'' information about~$M$,
which can only be seen from the outside$\ldots\,$.
The $a$ which we construct will be referred to as
a ``generic'' set relative to~$M$.
The idea is that all the properties of $a$ must be ``forced''
to hold merely on the basis that $a$ behaves like a ``generic''
set in~$M$.
This concept of deciding when a statement about~$a$
is ``forced'' to hold is the key point of the construction.
\end{quotation}
Cohen then proceeds to explain the forcing concept,
but at this point we will diverge from Cohen's account
and pursue instead the concept of a \textit{Boolean-valued model}
of~ZFC.
This approach was developed by
Scott, Solovay, and Vop\v enka starting in 1965,
and in my opinion is the most intuitive way to proceed at this juncture.
We will return to Cohen's approach later.

\section{Boolean-valued models}

To recap, we have reached the point where we see that if we want to
construct a model of $\neg$CH, it would be nice to have a method
of starting with an arbitrary standard transitive model~$M$ of~ZFC,
and building a new structure
by adjoining some subsets that are missing from~$M$.
We explain next how this can be done,
but instead of giving the construction right away,
we will work our way up to it gradually.

\subsection{Motivational discussion}

Inspired by Cohen's suggestion, we begin by
considering \textit{all possible statements} that might be true of
our new structure,
and then deciding which ones we want to hold and which one we do not
want to hold.

To make the concept of ``all possible statements'' precise,
we must introduce the concept of a \textit{formal language.}
Let $\mathfrak S$ denote the set of all sentences in the
\textit{first-order language of set theory,} i.e., all sentences
built out of ``atomic'' statements such as $x = y$ and $x R y$
(where $x$ and~$y$ are \textit{constant symbols} that each represent
some fixed element of the domain)
using the Boolean connectives OR, AND, and NOT
and the quantifiers $\exists$ and $\forall$.
The axioms of ZFC can all be expressed in this formal language,
as can any theorems (or non-theorems, for that matter) of ZFC.
For example, ``if $A$ then $B$'' (written $A \to B$) can be
expressed as $(\mathrm{NOT}\,A) \OR B$,
``$A$ iff $B$'' (written $A\leftrightarrow B$)
can be expressed as $(A\to B) \AND (B \to A)$,
``$x$ is a subset of~$y$'' (written $x \subseteq y$)
can be expressed as $\forall z\, ((z R x) \to (z R y))$,
and the powerset axiom can be expressed as
\begin{equation*}
\forall x \, \exists y \, \forall z \,
((z \subseteq x) \leftrightarrow (z R y)).
\end{equation*}

An important observation is that
when choosing which sentences in~$\mathfrak S$ we want to hold
in our new structure,
we are subject to certain constraints.
For example, if the sentences $\phi$ and $\psi$ hold,
then the sentence $\phi \AND \psi$ must also hold.
A natural way to track these constraints is by means of
a \textit{Boolean algebra.}
The most familiar example of a Boolean algebra is the family~$2^S$ of
all subsets of a given set~$S$, partially ordered by inclusion.
More generally, a Boolean algebra is any partially ordered set
with a minimum element~$\mathbf 0$ and a maximum element~$\mathbf 1$,
in which any two elements $x$ and~$y$ have a least upper bound $x\vee y$
and a greatest lower bound $x\wedge y$
(in the example of $2^S$, $\vee$ is set union
and $\wedge$ is set intersection),
where $\vee$ and $\wedge$ distribute over each other
(i.e., $x\vee(y\wedge z) = (x\vee y) \wedge (x \vee z)$
and $x\wedge(y\vee z) = (x\wedge y) \vee (x\wedge z)$)
and every element $x$ has a \textit{complement,}
i.e., an element $x^*$ such that
$x\vee x^* = \mathbf 1$ and $x\wedge x^* = \mathbf 0$.

There is a natural correspondence between the concepts
$\mathbf 0$, $\mathbf 1$, $\vee$, $\wedge$, and~$*$ in a Boolean algebra
and the concepts of falsehood, truth, OR, AND, and NOT in logic.
This observation suggests that
if we know that we want certain statements of~$\mathfrak S$
to hold in our new structure but are unsure of others,
then we can try to record our
state of partial knowledge by picking a suitable Boolean algebra~$\bB$,
and mapping every sentence $\phi \in \mathfrak S$
to some element of~$\bB$ that we denote by $\bval\phi$.
If $\phi$ is ``definitely true'' then we set
$\bval\phi = \mathbf 1$
and if $\phi$ is ``definitely false'' then we set
$\bval\phi = \mathbf 0$;
otherwise, $\bval\phi$ takes on some intermediate value
between $\mathbf 0$ and~$\mathbf 1$.
In a sense, we are developing a kind of ``multi-valued logic''
or ``fuzzy logic''\footnote{We use scare quotes as these
terms, and the term ``fuzzy set'' that we use later,
have meanings in the literature that are rather different from the
ideas that we are trying to convey here.}
in which some statements are neither true nor false
but lie somewhere in between.

It is clear that the mapping $\phi \mapsto \bval\phi$
should satisfy the conditions
\begin{align}
\label{eq:bvmvee}
\bval{\phi \OR \psi} &= \bval\phi \vee \bval\psi \\
\label{eq:bvmwedge}
\bval{\phi \AND \psi} &= \bval\phi \wedge \bval\psi \\
\label{eq:bvmsim}
\bval{\mathrm{NOT}\, \phi} &= (\bval\phi)^*
\end{align}
What about atomic expressions such as $\bval{x = y}$
and $\bval{x R y}$?
Again, if we definitely want certain equalities or membership statements
to hold but want to postpone judgment on others,
then we are led to the idea of tracking these statements
using a structure consisting of ``fuzzy sets.''
To make this precise, let us first observe that
an ordinary set may be identified with a function whose range is
the trivial Boolean algebra with just two elements $\mathbf 0$ and~$\mathbf 1$,
and that sends the members of the set to~$\mathbf 1$
and the non-members to~$\mathbf 0$.
Generalizing, if $\bB$ is an arbitrary Boolean algebra,
then a ``fuzzy set'' should take a set of
``potential members,'' which should themselves be fuzzy sets,
and assign each potential member~$y$ a value in~$\bB$
corresponding to the ``degree'' to which $y$ is a member of~$x$.
More precisely, we define a a \textit{$\bB$-valued set}
to be a function from a set of $\bB$-valued sets to~$\bB$.
(Defining $\bB$-valued sets in terms of $\bB$-valued sets
might appear circular, but the solution
is to note that the empty set is a $\bB$-valued set;
we can then build up other $\bB$-valued sets inductively.)

\subsection{Construction of $M^\bB$}

We are now in a position to describe more precisely our plan 
for constructing a new model of ZFC from a given model~$M$.
We pick a suitable Boolean algebra~$\bB$, and we let
$M^\bB$ be the set of all $\bB$-valued sets in~$M$.
The set $\mathfrak S$ should have one constant symbol
for each element of~$M^\bB$.
We define a a map $\phi \mapsto \bval\phi$ from
$\mathfrak S$ to~$\bB$, which should obey equations
such as (\ref{eq:bvmvee})--(\ref{eq:bvmsim})
and should send the axioms of ZFC to~$\mathbf 1$.
The structure $M^\bB$ will be a so-called \textit{Boolean-valued model}
of~ZFC; it will not actually be a model of~ZFC,
because it will consist of ``fuzzy sets'' and not sets,
and if you pick an arbitrary $\phi\in \mathfrak S$
and ask whether it holds in~$M^\bB$,
then the answer will often be neither ``yes'' nor ``no''
but some element of~$\bB$
(whereas if $N$ is an actual model of~ZFC
then either $N$ satisfies $\phi$ or it doesn't).
On the other hand,
$M^\bB$ will satisfy ZFC,
in the sense that $\bval\phi = \mathbf 1$ for every $\phi$ in~ZFC.
To turn $M^\bB$ into an actual model of ZFC with desired properties,
we will take a suitable quotient of~$M^\bB$ that eliminates the fuzziness.

We have already started to describe $M^\bB$ and the map
$\bval\cdot$, but we are not done.
For example, we need to deal with
expressions involving the quantifiers $\exists$ and~$\forall$.
These may not appear to have a direct counterpart in the
formalism of Boolean algebras, but notice that
another way to say that there exists $x$ with a certain property
is to say that either $a$ has the property or $b$ has the property
or $c$ has the property or$\ldots\,$, where we enumerate all the
entities in the universe one by one.
This observation leads us to the definition
\begin{equation}
\label{eq:bvmexists}
\bval{\exists x \, \phi(x)} = \bigvee_{a \in M^\bB} \bval{\phi(a)}
\end{equation}
Now there is a potential problem with (\ref{eq:bvmexists}):
In an arbitrary Boolean algebra,
an \textit{infinite} subset of elements may not have a least upper bound,
so the right-hand side of~(\ref{eq:bvmexists}) may not be defined.
We solve this problem by fiat: First we define a
\textit{complete Boolean algebra} to be a Boolean algebra in
which arbitrary subsets of elements have a least upper bound
and a greatest lower bound.
We then require that $\bB$ be a complete Boolean algebra;
then (\ref{eq:bvmexists}) makes perfect sense, as does the equation
\begin{equation}
\label{eq:bvmforall}
\bval{\forall x \, \phi(x)} = \bigwedge_{a \in M^\bB} \bval{\phi(a)}
\end{equation}

Equations (\ref{eq:bvmexists}) and (\ref{eq:bvmforall})
take care of $\exists$ and $\forall$,
but we have still not defined
\hbox{$\bval{x R y}$}
or
\hbox{$\bval{x = y}$},
or ensured that $M^\bB$ satisfies ZFC.
The definitions of $\bval{x = y}$
and $\bval{x R y}$ are surprisingly delicate;
there are many plausible attempts that fail for subtle reasons.
The impatient reader can safely skim the details in the next paragraph
and just accept the final equations
(\ref{eq:bvmmember})--(\ref{eq:bvmequal}).

We follow the treatment on pages 22--23 in Bell's book~\cite{bell},
which motivates the definitions of $\bval{x = y}$
and $\bval{x R y}$ by listing several equations
that one would like to hold and inferring what the definitions
``must'' be.
First, we want the axiom of extensionality to hold in~$M^\bB$;
this suggests the equation
\begin{equation*}
\bval{x = y} = 
\bval{ (\forall w\, (wRx \to wRy)) \AND
        (\forall w\, (wRy \to wRx)) }.
\end{equation*}
Another plausible equation is
\begin{equation*}
\bval{x R y} = \bval{ \exists w \, ((w R y) \AND (w = x)) }.
\end{equation*}
It is also plausible that the
expression $\bval{\exists w \, ((w R y) \AND \phi(w))}$
should depend only on the values of $\bval{\phi(w)}$
for those $w$ that are actually in the domain of~$y$
(recall that $y$, being a $\bB$-valued set,
is a function from a set of $\bB$-valued sets to~$\bB$,
and thus has a domain $\mathrm{dom}(y)$).
Also, the value of $\bval{w R y}$
should be closely related to the value of~$y(w)$.
We are thus led to the equations
\begin{align*}
\bval{ \exists w \, (w R y \AND \phi(w)) } &= 
\bigvee_{w\in\mathrm{dom}(y)}
  \bigl(y(w) \wedge \bval{ \phi(w) } \bigr) \\
\bval{ \forall w \, (w R y \to \phi(w)) } &= 
\bigwedge_{w\in\mathrm{dom}(y)}
  \bigl(y(w) \Rightarrow \bval{ \phi(w) }\bigr)
\end{align*}
where $x\Rightarrow y$ is another way of writing $x^* \vee y$.
All these equations drive us to the definitions
\begin{align}
\label{eq:bvmmember}
\bval{x R y} &= 
\bigvee_{w\in\mathrm{dom}(y)}
  \bigl(y(w) \wedge \bval{x = w} \bigr) \\
\label{eq:bvmequal}
\bval{x = y} &= 
\bigwedge_{w\in\mathrm{dom}(x)}
  \bigl(x(w) \Rightarrow \bval{ w R y } \bigr) \wedge
\bigwedge_{w\in\mathrm{dom}(y)}
  \bigl(y(w) \Rightarrow \bval{ w R x } \bigr)
\end{align}
The definitions (\ref{eq:bvmmember}) and (\ref{eq:bvmequal})
again appear circular,
because they define 
$\bval{ x = y}$ and $\bval{ x R y }$
in terms of each other,
but again (\ref{eq:bvmmember}) and (\ref{eq:bvmequal})
should be read as a joint inductive definition.

One final remark is needed regarding the definition of~$M^\bB$.
So far we have not imposed any constraints on~$\bB$
other than that it be a complete Boolean algebra.
But without some such constraints, there is no guarantee
that $M^\bB$ will satisfy ZFC.
For example, let us see what happens with the powerset axiom.
Given $x$ in $M^\bB$, it is natural to construct the
powerset~$y$ of~$x$ in~$M^\bB$ by letting
\begin{equation*}
\label{eq:bvmpowerset}
\mathrm{dom}(y) = \bB^{\mathrm{dom}(x)},
\end{equation*}
i.e., the ``potential members'' of~$y$ should be precisely the
maps from $\mathrm{dom}(x)$ to~$\bB$.
Moreover, for each $w\in \mathrm{dom}(y)$, the value of $y(w)$
should be $\bval{w\subseteq x}$.
The catch is that if $\bB$ is not in~$M$,
then maps from $\mathrm{dom}(x)$ to~$\bB$
may not be $\bB$-valued sets in~$M$.
The simplest way out of this difficulty is to
require that $\bB$ be in~$M$,
and we shall indeed require this.\footnote{While
choosing $\bB$ to be in~$M$ suffices to make everything work,
it is not strictly necessary.
\textit{Class forcing} involves certain carefully constructed
Boolean algebras $\bB$ that are not in~$M$.
However, this is an advanced topic that is
not needed for proving the independence of~CH.}
Once we impose this condition,
we can weaken
the requirement that $\bB$ be a complete Boolean algebra to
the requirement that $\bB$ be a complete Boolean algebra
\textit{in~$M$,} meaning that
infinite least upper bounds and greatest lower bounds
\textit{over subsets of~$B$ that are in~$M$} are guaranteed to exist,
but not necessarily in general.
(``Complete,'' being related to taking powersets, is not absolute.)
Examination of the definitions
of $M^\bB$ and $\bval\cdot$
reveals that $\bB$ only needs to be a complete Boolean algebra in~$M$,
and it turns out that
this increased flexibility in the choice of~$\bB$ is very important.

We are now done with the definition of the Boolean-valued model~$M^\bB$.
To summarize, we pick a Boolean algebra~$\bB$ in~$M$
that is complete in~$M$,
let $M^\bB$ be the set of all $\bB$-valued sets in~$M$,
and define $\bval\cdot$ using equations
(\ref{eq:bvmvee})--(\ref{eq:bvmequal}).

At this point, one needs to perform a long verification
that $M^\bB$ satisfies ZFC,
and that the rules of logical inference behave as expected in~$M^\bB$
(so that, for example, if $\bval\phi = \mathbf 1$ and
$\psi$ is a logical consequence of~$\phi$
then $\bval\psi = \mathbf 1$).
We omit these details because they are covered well in Bell's
book~\cite{bell}.
Usually, as in the case of the powerset axiom above,
it is not too hard to guess how to construct the object whose
existence is asserted by the ZFC axiom,
using the fact that $M$ satisfies ZFC,
although in some cases, completing the argument in detail can be tricky.

\subsection{Modding out by an ultrafilter}

As we stated above, the way to convert our Boolean-valued model~$M^\bB$
to an actual model of~ZFC is to take a suitable quotient.
That is, we need to pick out precisely the statements that
are true in our new model.
To do this, we choose a subset~$U$ of~$\bB$
that contains $\bval\phi$
for every statement $\phi$ that holds in the new model of~ZFC.
The set $U$, being a ``truth definition'' for our new model,
has to have certain properties;
for example, since for every~$\phi$, either $\phi$ or $\mathrm{NOT}\,\phi$
must hold in the new model,
it follows that for all~$x$ in~$\bB$,
$U$ must contain either $x$ or~$x^*$.
Similarly, thinking of membership in~$U$ as representing ``truth,''
we see that $U$ should have the following properties:
\begin{enumerate}
\item $1\in U$;
\item $0\notin U$;
\item if $x\in U$ and $y\in U$ then $x\wedge y \in U$;
\item if $x\in U$ and $x \le y$ (i.e., $x \wedge y = x$)
then $y \in U$;
\item for all $x$ in $\bB$, either $x\in U$ or $x^* \in U$.
\end{enumerate}
A subset $U$ of a Boolean algebra having the above properties
is called an \textit{ultrafilter.}

Given any ultrafilter~$U$ in~$\bB$
($U$ does not have to be in~$M$),
we define the quotient $\MBU$ as follows.
The elements of $\MBU$ are equivalence classes of elements of~$M^\bB$
under the equivalence relation
\begin{equation*}
x \sim_U y \quad \mbox{iff} \quad \bval{ x = y } \in U.
\end{equation*}
If we write $x^U$ for the equivalence class of~$x$, then
the binary relation of~$\MBU$---which we shall denote
by the symbol~$\in_U$---is defined by
\begin{equation*}
x^U \in_U y^U \quad \mbox{iff} \quad \bval{ x R y } \in U.
\end{equation*}
It is now fairly straightforward to verify that $\MBU$ is a model of~ZFC;
the hard work has already been done in verifying that $M^\bB$ satisfies~ZFC.

\section{Generic ultrafilters and the conclusion of the proof sketch}

At this point we have a powerful theorem in hand.
We can take any model~$M$,
any complete Boolean algebra~$\bB$ in~$M$,
and any ultrafilter $U$ of~$M$, and form a new model~$\MBU$ of~ZFC.
We can now experiment with various choices of $M$, $\bB$, and~$U$
to construct all kinds of models of~ZFC with various properties.

So let us revisit our plan (in Section~\ref{sec:CH}) of starting with
a standard transitive model and inserting some missing subsets
to obtain a larger standard transitive model.
If we try to use our newly constructed machinery to carry out this plan,
then we soon find that $\MBU$ need not, in general,
be (isomorphic to) a standard transitive model of~ZFC,
even if $M$ is.  Some extra conditions need to be imposed.

Cohen's insight---perhaps his most important and ingenious one---is
that in many cases, including the case of~CH,
the right thing to do is to require that $U$ be \textit{generic.}
The term ``generic'' can be defined more generally
in the context of an arbitrary partially ordered set~$P$.
First define a subset $D$ of~$P$ to be \textit{dense}
if for all $p$ in~$P$, there exists $q$ in~$D$ such that $q\le p$.
Then a subset of~$P$ is \textit{generic} if it intersects every dense subset.
In our current setting,
the partially ordered set~$P$ is $\bB \backslash \{{\mathbf 0}\}$,
and the crucial condition on~$U$ is that it be \textit{$M$-generic}
(or \textit{generic over~$M$}),
meaning that $U$ intersects
every dense subset~$D \subseteq (\bB \backslash \{{\mathbf 0}\})$
that is a member of~$M$.

If $U$ is $M$-generic, then $\MBU$ has many nice properties;
it is (isomorphic to) a standard transitive model of~ZFC,
and equally importantly, it contains~$U$.
In fact, if $U$ is $M$-generic,
then $\MBU$ is the smallest standard transitive model of~ZFC
that contains both $M$ and~$U$.
For this reason, when $U$ is $M$-generic, one typically writes
$M[U]$ instead of~$\MBU$.
We have realized the dream of adjoining a new subset of~$M$
to obtain a larger model
(remember that $U$ is a subset of~$\bB$ and
we have required $\bB$ to be in~$M$).

It is, of course, not clear that $M$-generic ultrafilters exist in general.
However, if $M$ is countable, then it turns out to be easy to prove
the existence of $M$-generic ultrafilters;
essentially, one just lists the dense sets and hits them one by one.
If $M$ is uncountable then
the Boolean-valued model machinery still works fine,
but $M$-generic ultrafilters may not exist.\footnote{This limitation
of uncountable models is not a big issue in
practice, because typically the L\"owenheim-Skolem theorem
allows us to replace an uncountable model with a countable surrogate.}
Fortunately for us,
the idea sketched at the beginning of Section~\ref{sec:CH}
relies on $M$ being countable anyway.

Let us now return to that idea and complete the proof sketch.
Start with a countable standard transitive model~$M$ of~ZFC.
If $M$ does not already satisfy~$\neg$CH,
then let $P$ be the partially ordered set of all
\textit{finite partial functions} from $\aleph_2^M\times \aleph_0$ into~$2$,
partially ordered by \textit{reverse} inclusion.
(A finite partial function is a finite set of ordered pairs
whose first coordinate is in the domain and whose second coordinate
is in the range, with the property that no two ordered pairs
have the same first element.)
There is a standard method, which we shall not go into here,
of \textit{completing} an arbitrary partially ordered set to a
complete Boolean algebra;
we take the completion of~$P$ in~$M$ to be our Boolean algebra~$\bB$.
Now take an $M$-generic ultrafilter~$U$,
which exists because $M$ is countable.
If we blur the distinction between $P$ and its completion~$\bB$ for a moment,
then we claim that $F := \bigcup U$ is a partial function from
$\aleph_2^M\times \aleph_0$ to~$2$.
To check this, we just need to check that any two elements $x$ and~$y$ of~$U$
are consistent with each other where they are both defined,
but this is easy: Since $U$ is an ultrafilter,
$x$ and~$y$ have a common lower bound $z$ in~$U$,
and both $x$ and~$y$ are consistent with~$z$.
Moreover, $F$ is a \textit{total} function;
this is because $U$ is generic,
and the finite partial functions that are defined at a specified point
in the domain form a dense set
(we can extend any partial function by defining it at that
point if it is not defined already).
Also, the sequence of functions from $\aleph_0$ to~$2$ encoded by~$F$
are pairwise distinct;
again this is because $U$ is generic,
and the condition of being pairwise distinct is a dense condition.
The axioms of ZFC ensure that $F \in M[U]$,
so $M[U]$ gives us the desired model of~$\neg$CH.

There is one important point that we have swept under the rug in the above
proof sketch.
The set $\aleph_2^M$ is still hanging around in~$M[U]$,
but it is conceivable that $\aleph_2^M$
may no longer play the role of~$\aleph_2$ in~$M[U]$;
i.e., it may be that $\aleph_2^M \ne \aleph_2^{M[U]}$.
Cardinalities are not absolute,
and so \textit{cardinal collapse} can occur,
i.e., the object that plays the role of a particular cardinal number in~$M$
may not play that same role in an extension of~$M$.
In fact, cardinal collapse does not occur in this particular case
but this fact must be checked.\footnote{The fact that cardinals
do not collapse here can be traced to the fact that the
Boolean algebra in question satisfies a combinatorial condition
called the \textit{countable chain condition} in~$M$.}
We omit the details, since they are covered thoroughly in textbooks.

\section{But wait---what about forcing?}

The reader may be surprised---justifiably so---that we have come to the
end of our proof sketch without ever precisely defining \textit{forcing.}
Does ``forcing'' not have a precise technical meaning?

Indeed, it does.
In Cohen's original approach,
he asked the following fundamental question.
Suppose that we want to adjoin a ``generic'' set~$U$ to~$M$.
What properties of the new model will be ``forced'' to hold
if we have only partial information about~$U$, namely we
know that some element $p$ of~$M$ is in~$U$?

If we are armed with the machinery of Boolean-valued models,
then we can answer Cohen's question.
Let us informally say
that \textit{$p$ forces~$\phi$} (written $p\forces \phi$) if
for every $M$-generic ultrafilter~$U$, $\phi$ must hold in~$M[U]$
whenever $p\in U$.
Note that $U$ plays two roles simultaneously;
it \textit{is} the generic set that we are adjoining to~$M$,
and it also picks out the true statements in~$M[U]$.
By the definition of an ultrafilter,
we see that if $p\le \bval\phi$,
then $\phi$ must be true in~$M[U]$ if $p\in U$.
Therefore we can give the following formal definition of
``$p\forces \phi$'':
\begin{equation}
\label{eq:forcing}
p\forces \phi \quad\mbox{iff}\quad p \le \bval\phi.
\end{equation}
The simplicity of equation~(\ref{eq:forcing}) explains why
our proof sketch did not need to refer to forcing explicitly.
Forcing is actually implicit in the proof, but
since $\forces$ has such a simple definition in terms
of~$\bval\cdot$,
it is possible in principle to produce a proof of Cohen's result
without explicitly using the symbol $\forces$ at all,
referring only to $\bval\cdot$
and Boolean algebra operations.

Of course,
Cohen did not have the machinery of Boolean-valued models available.
What he did was to figure out what properties
the expression $p\forces\phi$ ought to have,
given that one is trying to capture the notion of the logical implications
of knowing that $p$ is a member of our new ``generic'' set.
For example, one should have
$p\forces (\phi\AND\psi)$ iff $p\forces\phi$ and $p\forces\psi$,
by the following reasoning:
If we know that membership of~$p$ in~$U$
forces $\phi$ to hold and it also forces $\psi$ to hold,
then membership of $p$ in~$U$ must also force $\phi\AND\psi$ to hold.

By similar but more complicated reasoning,
Cohen devised a list of rules analogous to
(\ref{eq:bvmvee})--(\ref{eq:bvmequal}) that he used to define
$p\forces\phi$ for any statement~$\phi$ in~$\mathfrak S$.
In this way, he built all the necessary machinery
on the basis of the forcing relation,
without ever having to introduce Boolean algebras.

Thus there are (at least) two different ways to approach this subject,
depending on whether $\forces$ or $\bval\cdot$
is taken to be the fundamental concept.
For many applications, these two approaches ultimately
amount to almost the same thing,
since we can use equation~(\ref{eq:forcing}) to pass between them.
In this paper I have chosen the approach using
Boolean-valued models because I feel that
the introduction of ``fuzzy sets'' and the Boolean algebra~$\bB$
are relatively easy to motivate.
In Cohen's approach one still needs to introduce at some point
some ``fuzzy sets'' (called \textit{names} or \textit{labels})
and a partial order,
and these seem (to me at least) to be pulled out of a hat.
Also, the definitions (\ref{eq:bvmvee})--(\ref{eq:bvmequal})
are somewhat simpler than the corresponding definitions for~$\forces$.

On the other hand, even when one works with Boolean-valued models,
Cohen's intuition about generic sets~$U$,
and what is forced to be true if we know that $p\in U$,
is often extremely helpful.
For example, recall from our proof sketch that the constructed
functions from $\aleph_0$ to~$2$ were pairwise distinct.
In geometry, ``generically'' chosen functions will not be equal;
distinctness is a dense condition.
Cohen's intuition thus leads us to the (correct) expectation
that our generically chosen functions
will also be distinct, because no finite $p\in U$
can force two of them to be equal.
This kind of reasoning is invaluable in more complicated applications.

\section{Final remarks}

We should mention that the Boolean-valued-model approach has
some disadvantages.
For example, set theorists sometimes find the need to work
with models of axioms that do not include the powerset axiom,
and then the Boolean-valued model approach does not work,
because ``complete'' does not really make sense in such contexts.
Also, Cohen's original approach
allows one to work directly with an arbitrary partially ordered set~$P$
that is not necessarily a Boolean algebra,
and a \textit{generic filter} rather than a generic ultrafilter.
(A subset~$F$ of~$P$ is a \textit{filter} if $p\in F$ and $p\le q$
implies $q\in F$, and every $p$ and~$q$ have a common lower bound in~$F$.)
In our proof sketch we have already caught a whiff of the fact that
in many cases, there is some partially ordered set~$P$ lying around that
captures the combinatorics of what is really going on,
and having to complete~$P$ to a Boolean algebra is a technical nuisance;
it is much more convenient to work with $P$ directly.
If the reader prefers this approach, then
Kunen's book~\cite{kunen} would be my recommended reference.
Note that Kunen helpfully supplements his treatment with
an abbreviated discussion of Boolean-valued models
and the relationship between the two different approaches.

In my opinion, the weakest part of the exposition in this paper
is the treatment of \textit{genericity,}
whose definition appears to come out of nowhere.
A~posteriori one can see that the definition works beautifully,
but how would one guess a~priori that the geometric concepts
of dense sets and generic sets would be so apropos in this context,
and come up with the right precise definitions?
Perhaps the answer is just that Cohen was a genius,
but perhaps there is a better approach
yet to be discovered that will make it all clear.

Let us conclude with some suggestions for further reading.
Easwaran \cite{easwaran} and Wolf \cite{wolf} give
very nice overviews of forcing written in the same spirit as
the present paper, giving details that are critical for understanding
but omitting messy technicalities.
Scott's paper \cite{scott} is
a classic exposition written for non-specialists,
and Cohen \cite{cohen2} gave a lecture late in his life
about how he discovered forcing.
The reader may also find it helpful to study the connections
that forcing has with topology \cite{bowen},
topos theory \cite{mm}, modal logic \cite{smullyanfitting},
arithmetic \cite{boolosjeffrey}, proof theory \cite{avigad} and
computational complexity \cite{krajicek}.
It may be that insights from these differing perspectives
can be synthesized to solve the open exposition problem of forcing.

\section{Acknowledgments}

This paper grew out of an article entitled ``Forcing for dummies''
that I posted to the USENET newsgroup \texttt{sci.math.research} in 2001.
However, in that article I did not employ the Boolean-valued model approach,
and hence the line of exposition was quite different.
Interested readers can easily find the earlier version on the web.

I thank
Costas Drossos,
Ali Enayat,
Martin Goldstern, and
Matt Yurkewych
for pointing me to references
\cite{scott},
\cite{bowen},
\cite{cohen2}, and
\cite{smullyanfitting}
respectively.

I am deeply indebted to Matthew Wiener and especially Andreas Blass
for their many prompt and patient answers to my dumb questions about forcing.
Thanks also to Fred Kochman, Miller Maley, Christoph Weiss, Steven Gubkin,
James Hirschhorn, and the referee
for comments on an earlier draft of this paper that have improved the
exposition and fixed some bugs.

\end{document}